2005, Vol. 15, No. 1A, 125–152

# STRUCTURE OF LARGE RANDOM HYPERGRAPHS

By R. W. R. Darling and J. R. Norris

*National Security Agency and University of Cambridge*

The theme of this paper is the derivation of analytic formulae for
certain large combinatorial structures. The formulae are obtained via
fluid limits of pure jump-type Markov processes, established under
simple conditions on the Laplace transforms of their Lévy kernels.
Furthermore, a related Gaussian approximation allows us to describe
the randomness which may persist in the limit when certain param-
eters take critical values. Our method is quite general, but is ap-
plied here to vertex identifiability in random hypergraphs. A vertex
$v$ is identifiable in $n$ steps if there is a hyperedge containing $v$ all
of whose other vertices are identifiable in fewer steps. We say that
a hyperedge is identifiable if every one of its vertices is identifiable.
Our analytic formulae describe the asymptotics of the number of
identifiable vertices and the number of identifiable hyperedges for a
Poisson($\beta$) random hypergraph $\Lambda$ on a set $V$ of $N$ vertices, in the
limit as $N \to \infty$. Here $\beta$ is a formal power series with nonnegative
coefficients $\beta_0, \beta_1, \ldots,$ and $(\Lambda(A))_{A \subseteq V}$ are independent Poisson ran-
dom variables such that $\Lambda(A)$, the number of hyperedges on $A$, has
mean $N\beta_j / \binom{N}{j}$ whenever $|A| = j$.

## 1. Introduction.

1.1. *Motivation.* We are interested in the evolution of certain statisti-
cally symmetric random structures, extended over a large finite set of points,
when points are progressively removed in a way which depends on the struc-
ture. The initial condition of the structure may allow few possibilities for
the removal of points; indeed, it may be that, once a small proportion of
points are removed, the process terminates. On the other hand, the removal
of points may cause the structure to ripen, eventually yielding a large pro-
portion of the initial points. Our analysis will enable us to demonstrate a

Received November 2002; revised May 2004.

graph.

This is an electronic reprint of the original article published by the
Institute of Mathematical Statistics in *The Annals of Applied Probability*,
2005, Vol. 15, No. 1A, 125–152. This reprint differs from the original in
pagination and typographic detail.

sharp transition between these two sorts of behavior as certain parameters pass through critical values.

Let us illustrate this phenomenon by a simple special case. Consider the complete graph on $N$ vertices and declare each vertex to be open with probability $p$, each edge to be open with probability $\alpha/N$. Suppose that we are allowed to select an open vertex, remove it, and declare open any other vertices sharing an open edge with the selected vertex. If we continue in this way until no open vertices remain, we eventually remove every vertex connected to an open vertex by open edges. We shall see that the proportion of vertices thus removed converges in probability as $N \to \infty$ and that the limit $z^*(p, \alpha)$ is the unique root in $[0, 1)$ of the equation

$$\alpha z + \log(1 - z) = \log(1 - p).$$

Thus, for small values of $p$, there is a dramatic change in behavior as $\alpha$ passes through 1. As $p \downarrow 0$, for $\alpha \leq 1$,

$$z^*(p, \alpha)/p \to 1/(1 - \alpha),$$

but for $\alpha > 1$,

$$z^*(0+, \alpha) > 0.$$

Of course, this is a reflection of well-known connectivity properties of random graphs, discovered by Erdős and Rényi [7], and discussed, for example, in [3].

The class of models considered in this paper is a natural generalization of some classical models of random graphs and hypergraphs, which may be further motivated as follows. Phase transitions in combinatorial problems constitute an area of active research among computer scientists. Many "hard" combinatorial problems can be cast as *satisfiability* problems, which seek to assign a truth value to each of a set of Boolean variables, such that a collection of logical conjunctions are simultaneously satisfied. Phase transitions for random satisfiability ("random k-SAT") problems have been studied by researchers at Microsoft [1, 2, 15] and IBM [4], but difficult questions remain unanswered. The random hypergraph model herein may be viewed as a simplification of the random satisfiability model: a vertex corresponds to a Boolean variable, and a hyperedge to the set of variables appearing in a specific logical conjunction, neglecting the truth or falsehood assigned to those variables. Under this simplification, definitive critical parameters are obtained which shed light on the random satisfiability model, and whose derivation may serve as a template for analysis of mixed satisfiability problems.

1.2. *Hypergraphs.* Let $V$ be a finite set of $N$ vertices. By a *hypergraph* on $V$ we mean any map

$$\Lambda : \mathcal{P}(V) \to \mathbb{Z}^+.$$

Here $\mathbb{Z}^+$ denotes the set of nonnegative integers. The reader may consult [6] for an overview of the theory of hypergraphs; however, the direction pursued here is largely independent of previous work. We emphasize that, in distinction to much of the combinatorial literature on hypergraphs, we allow the possibility that more than one edge is assigned to a given subset; thus we are considering multi-hypergraphs. Moreover, we do not insist that all hyperedges have the same number of vertices. Much of the literature is restricted to this *uniform* case. Our methods allow a significant broadening of the class of models for which asymptotic computations are feasible. Hyperedges over vertices are called *patches* (loops in [6]) and hyperedges over $\varnothing$ are called *debris*. The total number of hyperedges is

$$|\Lambda| = \sum_A \Lambda(A).$$

1.3. *Accessibility and identifiability.* Interest in large random graphs has often focused on the sizes of their connected components. If there is given also, as in the example above, a set of distinguished vertices $V_0$, then it is natural to seek to determine the proportion of all vertices connected to $V_0$.

In the more general context of hypergraphs there is more than one interesting counterpart of connectivity. Given a hypergraph $\Lambda$ on a set $V$, we say that a vertex $v$ is *accessible in 1 step* or, equivalently, *identifiable in 1 step* if $\Lambda(\{v\}) \geq 1$. We say, for $n = 2, 3, \ldots$, that a vertex is *accessible in $n$ steps* if it belongs to some subset $A$ with $\Lambda(A) \geq 1$, *some* other element of which is accessible in less than $n$ steps. A vertex is *accessible* if it is accessible in $n$ steps for some $n \geq 1$.

On the other hand, we say that a vertex is *identifiable in $n$ steps* if it belongs to some subset $A$ with $\Lambda(A) \geq 1$, *all* of whose other elements are identifiable in less than $n$ steps. A vertex is *identifiable* if it is identifiable in $n$ steps for some $n \geq 1$.

The notion of accessibility may be appropriate to some physical models similar to percolation, whereas identifiability is more relevant to knowledge-based structures. We shall examine only the notion of identifiability.

Given a hypergraph $\Lambda$ without patches and a distinguished vertex $v_0$, we say that a vertex $v$ is *accessible from $v_0$* if it is accessible in the hypergraph $\Lambda + 1_{\{\{v_0\}\}}$, that is, in the hypergraph obtained from $\Lambda$ by adding a single patch at $v_0$. Identifiability from $v_0$ is defined similarly. The set of vertices accessible from $v_0$ is the *component* of $v_0$, as studied in [6, 11, 12, 14]. The set of vertices identifiable from $v_0$ is the *domain* of $v_0$, as studied by Levin

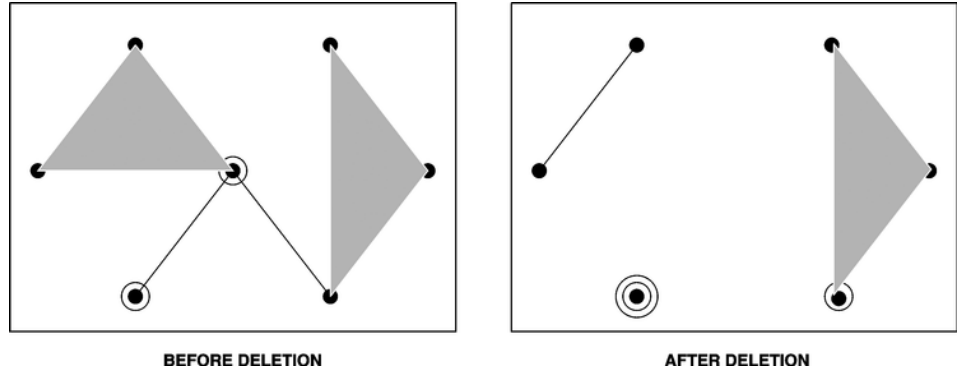

**BEFORE DELETION**                    **AFTER DELETION**

Fig. 1. *Example of a permitted collapse—deletion of one vertex.*

and the current authors [5]. We shall not consider further in this paper these vertex-based notions.

The process of identification is dual to the process leading to the 2-core of a graph or hypergraph, that is to say, the maximal subgraph in which every nonisolated vertex has degree at least 2. In the former process one removes vertices having a 1-hyperedge, in the latter one removes edges containing a vertex of degree 1. In this duality, nonidentifiable vertices correspond to the 2-core. Thus our results may be interpreted as giving the asymptotic size of the 2-core for a certain class of random hypergraphs.

1.4. *Hypergraph collapse.* It will be helpful to think of the identification of vertices as a progressive activity. Once a vertex is identified, it is *removed* or *deleted* from the vertex set, in a manner which is explained below. Thus, we shall consider an evolution of hypergraphs by the removal of vertices over which there is a patch. A hypergraph with no patches will therefore be *stable*. Given a hypergraph $\Lambda$ and a vertex $v$, we can arrive at a new hypergraph $\Lambda'$ by removing $v$ from each of the hyperedges of $\Lambda$. Thus

$$\Lambda'(A) = \begin{cases} \Lambda(A) + \Lambda(A \cup \{v\}), & \text{if } v \notin A, \\ 0, & \text{if } v \in A. \end{cases}$$

For example, in Figure 1, the patch on the central vertex is selected, and that vertex is removed; this causes a triangular face to collapse to an edge, and two edges incident to the vertex to collapse to patches on the vertices at the other ends. Note that this leaves two patches on the lower left vertex.

If $\Lambda(\{v\}) \geq 1$, then we say that $\Lambda'$ is obtained from $\Lambda$ by a (*permitted*) *collapse*. Starting from $\Lambda$, we can obtain, by a finite sequence of collapses, a stable hypergraph $\Lambda_\infty$. Denote by $V^*$ the set of vertices removed in passing from $\Lambda$ to $\Lambda_\infty$. The elements of $V^*$ are the identifiable vertices. We write $\Lambda^*$ for the *identifiable hypergraph*, given by

$$\Lambda^*(A) = \Lambda(A)\mathbb{1}_{A \subseteq V^*}.$$

We note that $V^*$, and hence $\Lambda^*$ and $\Lambda_\infty$, do not depend on the particular sequence of collapses chosen. For, if $v_1, v_2, \ldots$ and $v_1', v_2', \ldots$ are two such sequences, and if $v_n \neq v_k'$ for all $k$, then we can take $n$ minimal and find $k$ such that $\{v_1, \ldots, v_{n-1}\} \subseteq \{v_1', \ldots, v_k'\}$; then, with an obvious notation, $\Lambda_k'(\{v_n\}) \geq \Lambda_{n-1}(\{v_n\}) \geq 1$, so $v_n$ must, after all, appear in the terminating sequence $v_1', v_2', \ldots$, a contradiction. We note also that $V^*$ increases with $\Lambda$.

1.5. *Purpose of this paper.* The main question we shall address is to determine the asymptotic sizes of $V^*$ and $\Lambda^*$ for certain generic random hypergraphs, as the number of vertices becomes large. We note that, since the number of hyperedges is conserved in each collapse, all the identifiable hyperedges eventually turn to debris:

$$\Lambda_\infty(\varnothing) = |\Lambda^*|.$$

Note that $V^*$ depends only on $\min\{\Lambda, 1\}$. In the case where $\Lambda(A) = 0$ for $|A| \geq 3$, the hypergraph $\min\{\Lambda, 1\}$ may be considered as a graph on $V$ equipped with a number of distinguished vertices. Then $V^*$ is precisely the set of vertices connected in the graph to one of these distinguished vertices.

1.6. *Poisson random hypergraphs.* Let $(\Omega, \mathcal{F}, \mathbb{P})$ be a probability space. A *random hypergraph* on $V$ is a measurable map

$$\Lambda : \Omega \times \mathcal{P}(V) \to \mathbb{Z}^+.$$

An introduction to random hypergraphs may be found in [11], though we shall pursue rather different questions here. We shall consider a class of random hypergraphs whose distribution is determined by a sequence $\beta = (\beta_j : j \in \mathbb{Z}^+)$ of nonnegative parameters. Say that a random hypergraph $\Lambda$ on $V$ is Poisson$(\beta)$ if:

1. the random variables $\Lambda(A), A \subseteq V$, are independent;
2. the distribution of $\Lambda(A)$ depends only on $|A|$;
3. $\sum_{|A|=j} \Lambda(A) \sim \text{Poisson}(N\beta_j), j = 0, 1, \ldots, N$.

A consequence of these assumptions is that $\Lambda(A)$ has mean $N\beta_j / \binom{N}{j}$ whenever $|A| = j$. Note that, when $N$ is large, for $j \geq 2$, only a small fraction of the subsets of size $j$ have any hyperedges, and those that do usually have just one. Also the ratio of $j$-edges to vertices tends to $\beta_j$. Our assumption of Poisson distributions is a convenient exact framework reflecting behavior which holds asymptotically as $N \to \infty$ under more generic conditions.

1.7. *Generating function.* A key role is played by the power series

$$(1) \qquad \qquad \beta(t) = \sum_{j \geq 0} \beta_j t^j$$

and by the derived series

$$\beta'(t) = \sum_{j \geq 1} j \beta_j t^{j-1},$$

$$\beta''(t) = \sum_{j \geq 2} j(j-1) \beta_j t^{j-2}.$$

Let $\beta$ have radius of convergence $R$. The function $\beta'(t) + \log(1-t)$ may have zeros in $[0, 1)$ but these can accumulate only at 1. Set

(2)                 $z^* = \inf\{t \in [0, 1) : \beta'(t) + \log(1-t) < 0\} \wedge 1$

and denote by $\zeta$ the set of zeros of $\beta'(t) + \log(1-t)$ in $[0, z^*)$. Note that if $\beta$ is a polynomial, or indeed if $R > 1$, then $z^* < 1$. Also, the generic and simplest case is where $\zeta$ is empty.

**2. Results.** We state our principal result first in the generic case.

2.1. *Hypergraph collapse—generic case.*

THEOREM 2.1. *Assume that $z^* < 1$ and $\zeta = \varnothing$. For $N \in \mathbb{N}$, let $V^N$ be a set of $N$ vertices and let $\Lambda^N$ be a Poisson$(\beta)$ hypergraph on $V^N$. Then, as $N \to \infty$, the numbers of identifiable vertices and identifiable hyperedges satisfy the following limits in probability:*

$$|V^{N*}|/N \to z^*, \qquad |\Lambda^{N*}|/N \to \beta(z^*) - (1 - z^*)\log(1 - z^*).$$

EXAMPLE 2.1. The random graph with distinguished vertices described in the Introduction corresponds to a Poisson$(\beta^N)$ hypergraph $\Lambda^N$, where

$$1 - e^{-\beta_1^N} = p, \qquad 1 - e^{-2\beta_2^N/(N-1)} = \alpha/N$$

and $\beta_j^N = 0$ for $j \geq 3$. Note that $\beta_1^N = \beta_1$ and $\beta_2^N \to \beta_2$ as $N \to \infty$, where $\beta_1 = -\log(1-p)$ and $\beta_2 = \alpha/2$. Theorem 2.1 extends easily to cases where $\beta$ depends on $N$ in such a mild way: one just has to check that Lemma 6.1 remains valid and note that this is the only place that $\beta$ enters the calculations. We have $\beta(t) = -t\log(1-p) + t^2\alpha/2$ so

$$\beta'(t) + \log(1-t) = -\log(1-p) + t\alpha + \log(1-t).$$

Then $z^*$ is the unique $t \in [0, 1)$ such that

$$\alpha t + \log(1-t) = \log(1-p)$$

and $\zeta$ is empty, so $|V^{*N}|/N \to z^*$ in probability as $N \to \infty$, as stated above.

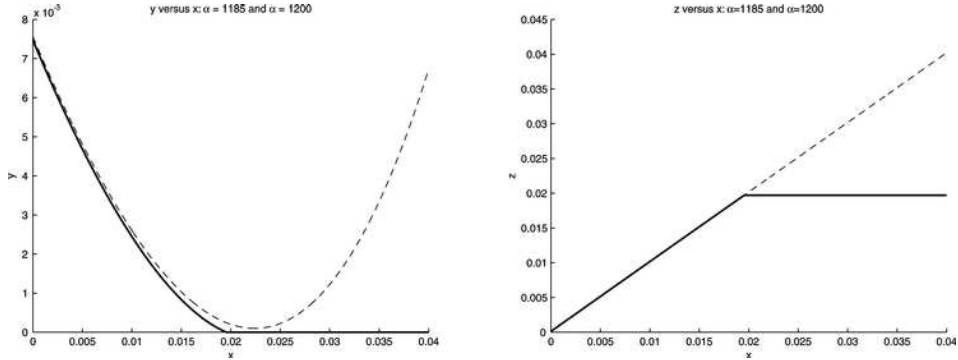

Fig. 2. *Changes in behavior near a critical parameter value.*

Example 2.2. To illustrate critical phenomena, let $\beta(t) = \alpha(0.1 + 0.9t)^7$. Let $x$, $y$ and $z$ refer to the re-scaled number of vertices eliminated, the number of patches and the amount of debris, respectively; here "re-scaled" means after division by the number of vertices. Plots of $y$ and $z$ versus $x$ are shown in Figure 2, for the choices $\alpha = 1185$ (solid) and $\alpha = 1200$ (dashed). In the case $\alpha = 1185$, $y$ hits zero when $x \approx 0.02$, and so $z$ remains stuck at about 0.02. A very small increase in $\alpha$, from 1185 to 1200, causes a dramatic change in the outcome: after narrowly avoiding extinction (Figure 2), the number of patches explodes (Figure 3) as $x$ increases toward 1.

Consider what the figures tell us about the supercritical case $\alpha = 1200$: during the first 4% of patch selections, there is rarely any other patch covering the same vertex as the one selected; Figure 3 shows that, during the last 10% of patch selections, an average of 5792 other patches cover the same vertex as the one selected. (Read the labels on the $x$-axes carefully: Figure 2 is a close-up of the leftmost 4% of the scale of Figure 3.)

2.2. *Hypergraph collapse—general case.* In order to describe an extension of Theorem 2.1 to the case where $\zeta$ is nonempty, we introduce the random variable

$$Z = \min\{z \in \zeta : W(z/(1-z)) < 0\} \wedge z^*,$$

where $(W_t)_{t \geq 0}$ is a Brownian motion.

Theorem 2.2. *Assume that $R \notin \zeta$. Then, for $V^{N*}$ and $\Lambda^{N*}$ as in Theorem 2.1, the following limits exist in distribution:*

$$|V^{N*}|/N \to Z, \qquad |\Lambda^{N*}|/N \to \beta(Z) - (1-Z)\log(1-Z).$$

In the case where $\zeta$ has only a single point $\zeta_0 < z^*$, then $Z$ is equal to $\zeta_0$ with probability $\frac{1}{2}$ and equal to $z^*$ with probability $\frac{1}{2}$. We do not know what happens when $R \in \zeta$. Proofs will be given in Section 6.

**3. Randomized collapse.** We introduce here a particular random rule for choosing the sequence of moves by which a hypergraph is collapsed, which has the desirable feature that certain key statistics of the evolving hypergraph behave as Markov chains. It is by analysis of the asymptotics of these Markov chains as $N \to \infty$ that we are able to prove our main results.

3.1. *Induced hypergraph.* Let $\Lambda$ be a Poisson($\beta$) hypergraph. For $S \subseteq V$ with $|S| = n$, let $\Lambda^S$ be the hypergraph obtained from $\Lambda_0$ by removing all vertices in $S$. Thus, for $A \subseteq V \setminus S$ with $|A| = j$,

$$\Lambda^S(A) = \sum_{B \supseteq A, B \setminus S = A} \Lambda(B) \sim P(\lambda_j(N, n)),$$

where the Poisson parameter $\lambda_j(N, n)$ is computed as follows: there are $\binom{n}{i}$ ways to choose $S \cap B$ such that $|B| = i + j$, and the Poisson parameter of $\Lambda_0(B)$ is $N\beta_{j+i}/\binom{N}{i+j}$, so

$$\lambda_j(N, n) = N \sum_{i=0}^{n} \beta_{j+i} \binom{n}{i} \Big/ \binom{N}{i+j}.$$

Moreover, the random variables $\Lambda^S(A)$, $A \subseteq V \setminus S$, are independent.

3.2. *Rule for randomized collapse.* Recall that the sequence of vertices chosen to collapse a hypergraph is unimportant, provided we keep going until there are no more patches. However, we shall use a specific randomized rule which turns out to admit a description in terms of a finite-dimensional Markov chain. This leads to a randomized process of collapsing hypergraphs $(\Lambda_n)_{n \geq 0}$. This will prove to be an effective means to compute the numbers of identifiable vertices and identifiable hyperedges for $\Lambda_0$.

The process $(\Lambda_n)_{n \geq 0}$, together with a sequence of sets $(S_n)_{n \geq 0}$ such that $\Lambda_n = \Lambda^{S_n}$, is constructed as follows. Let $S_0 = \varnothing$ and $\Lambda_0 = \Lambda$. Suppose that $S_n$ and $\Lambda_n$ have been defined. If there are no patches in $\Lambda_n$, then $S_{n+1} = S_n$ and $\Lambda_{n+1} = \Lambda_n$. If there are patches in $\Lambda_n$, select one uniformly at random and denote by $v_{n+1}$ the corresponding vertex; then set $S_{n+1} = S_n \cup \{v_{n+1}\}$ and $\Lambda_{n+1} = \Lambda^{S_{n+1}}$.

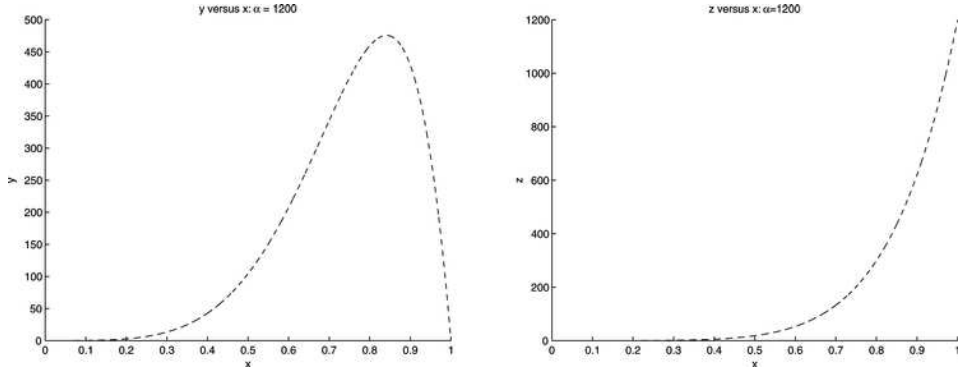

FIG. 3. *Patches and debris in the supercritical regime.*

**3.3.** *An embedded Markov chain.* Let $Y_n$ denote the number of patches and $Z_n$ the amount of debris in $\Lambda_n$. Then $Y_n = 0$ and $Z_n = |\Lambda^*|$ for $n \geq |V^*|$. Also $|V^*| = \inf\{n \geq 0 : Y_n = 0\}$. Let $W_{n+1}$ denote the number of other patches at time $n$ sharing the same vertex as the $(n+1)$st selected patch, and let $U_{n+1}$ denote the number of 2-edges at time $n$ containing the $(n+1)$st selected vertex $v_{n+1}$. Our analysis will rest on the observation that $(Y_n, Z_n)_{n \geq 0}$ is a Markov chain, where, conditional on $Y_n = m \geq 1$ and $Z_n = k$, we have

$$Y_{n+1} = Y_n - 1 - W_{n+1} + U_{n+1}, \qquad Z_{n+1} = Z_n + 1 + W_{n+1}$$

and where $W_{n+1} \sim B(m-1, 1/(N-n))$ and $U_{n+1} \sim P((N-n-1)\lambda_2(N,n))$ with $W_{n+1}$ and $U_{n+1}$ independent.

To see this, introduce the filtration

$$\mathcal{F}_n = \sigma(S_r, Y_r, Z_r : r = 0, 1, \ldots, n).$$

LEMMA 3.1. *Let*

$$p(\lambda | S, k, m) = \mathbb{P}\left[ \Lambda^S = \lambda \,\Big|\, \sum_v \Lambda^S(\{v\}) = m, \Lambda^S(\varnothing) = k \right].$$

*Then*

$$\mathbb{P}[\Lambda_n = \lambda | \mathcal{F}_n] = p(\lambda | S_n, Y_n, Z_n).$$

*Equivalently, for all $B \in \mathcal{F}_n$ so that $B \subset \{S_n = S, Y_n = m, Z_n = k\}$,*

$$\mathbb{P}[\Lambda_n = \lambda, B] = p(\lambda | S, m, k)\mathbb{P}[B].$$

The claimed Markov structure for $(Y_n, Z_n)_{n \geq 0}$ follows easily.

PROOF OF LEMMA 3.1. The identity is obvious for $n = 0$. Suppose it holds for $n$. Let $B \subset \{S_n = S, Y_n = m, Z_n = k\}$. Take $x \in V \setminus S$, $m' \geq 1$ and $k' > k$. Set $S' = S \cup \{x\}$ and $B' = \{S_{n+1} = S', Y_{n+1} = m', Z_{n+1} = k'\} \cap B$. It will suffice to show, for all hypergraphs $\lambda'$ having $m'$ patches and amount of debris $k'$, that

$$\mathbb{P}[\Lambda_{n+1} = \lambda', B'] \propto p(\lambda' | S', m', k'),$$

where $\propto$ denotes equality up to a constant independent of $\lambda'$. But

$$\mathbb{P}[\Lambda_{n+1} = \lambda', B'] = \sum_\lambda \frac{k' - k}{m}\mathbb{P}[\Lambda_n = \lambda, B] \propto \sum_\lambda p(\lambda | S, m, k),$$

where the sum is over all hypergraphs $\lambda$ which collapse to $\lambda'$ on removing the vertex $x$. Let $Y^S, Z^S$ denote the number of patches, amount of debris in $\Lambda^S$, respectively. Since $Y^S$ and $Z^S$ are conditionally independent of $\Lambda^{S'}$ given $Y^{S'}$ and $Z^{S'}$,

$$\sum_\lambda p(\lambda|S, m, k) = \mathbb{P}(\Lambda^{S'} = \lambda'|Y^S = m, Z^S = k) \propto p(\lambda'|S', m', k'),$$

as desired.  □

## 4. Exponential martingales for jump processes.

We recall here some standard notions for pure jump Markov processes in $\mathbb{R}^d$ and their associated martingales. These will be used to study the fluid limit of a sequence of such jump processes in Section 5.

4.1. *Laplace transforms.* Let $(X_t)_{t \geq 0}$ be a pure jump Markov process taking values in a subset $I$ of $\mathbb{R}^d$, with Lévy kernel $K$. Consider the Laplace transform

$$m(x, \theta) = \int_{\mathbb{R}^d} e^{\langle \theta, y \rangle} K(x, dy), \qquad \theta \in (\mathbb{R}^d)^*,$$

and assume that, for some $\eta_0 > 0$,

$$(3) \qquad\qquad \sup_{x \in I} \sup_{|\theta| \leq \eta_0} m(x, \theta) \leq C < \infty.$$

The distribution of the time $T$ and displacement $\Delta X_T$ of the first jump of $(X_t)_{t \geq 0}$ is given by

$$\mathbb{P}(T \in dt, \Delta X_T \in dy | T > t, X_0 = x) = K(x, dy)\, dt.$$

Introduce random measures $\mu$ and $\nu$ on $(0, \infty) \times \mathbb{R}^d$, given by

$$\mu = \sum_{\Delta X_t \neq 0} \varepsilon_{(t, \Delta X_t)},$$

$$\nu(dt, dy) = K(X_{t-}, dy)\, dt,$$

where $\varepsilon_{(t, y)}$ denotes the unit mass at $(t, y)$; $\nu$ is thus the compensator of the random measure $\mu$, in the sense of ([10], page 422).

4.2. *Martingales associated with jump processes.* The fact that $\nu$ is a compensator implies that, for any previsible process $a : \Omega \times (0, \infty) \times \mathbb{R}^d \to \mathbb{R}$ satisfying

$$\mathbb{E} \int_{\mathbb{R}^d} |a(s, y)| \nu(ds, dy) < \infty,$$

the following process is a martingale:

$$\int_0^t \int_{\mathbb{R}^d} a(s,y)(\mu - \nu)(ds, dy).$$

In particular, (3) allows us to take $a(s,y) = y$, which gives the martingale

$$M_t = \int_0^t \int_{\mathbb{R}^d} y(\mu - \nu)(ds, dy).$$

Fix $\eta \in (0, \eta_0)$. Then there exists $A < \infty$ such that

(4) $$|m''(x,\theta)| \leq A, \qquad x \in I, |\theta| \leq \eta,$$

where "$'$" denotes differentiation in $\theta$. Define for $\theta \in (\mathbb{R}^d)^*$

$$\phi(x,\theta) = \int_{\mathbb{R}^d} \{e^{\langle \theta, y \rangle} - 1 - \langle \theta, y \rangle\} K(x, dy).$$

Then $\phi \geq 0$ and, for $|\theta| \leq \eta$, by the second-order mean value theorem,

$$\phi(x,\theta) = \int_0^1 m''(x, r\theta)(\theta, \theta)(1 - r)\, dr$$

so

$$\phi(x,\theta) \leq \tfrac{1}{2} A|\theta|^2, \qquad x \in I, |\theta| \leq \eta.$$

Let $(\theta_t)_{t \geq 0}$ be a previsible process in $(\mathbb{R}^d)^*$ with $|\theta_t| \leq \eta$ for all $t$. Set

(5) $$Z_t = Z_t^\theta = \exp\Big\{ \int_0^t \langle \theta_s, dM_s \rangle - \int_0^t \phi(X_s, \theta_s)\, ds \Big\}.$$

Then $(Z_t)_{t \geq 0}$ is locally bounded, and by the Doléans formula ([10], page 440),

$$Z_t = 1 + \int_0^t \int_{\mathbb{R}^d} Z_{s-}(e^{\langle \theta_s, y \rangle} - 1)(\mu - \nu)(ds, dy).$$

Hence $(Z_t)_{t \geq 0}$ is a nonnegative local martingale, so $\mathbb{E}(Z_t) \leq 1$ for all $t$. Hence

$$\mathbb{E} \int_0^t \int_{\mathbb{R}^d} |Z_{s-}(e^{\langle \theta_s, y \rangle} - 1)| \nu(ds, dy)$$

$$\leq \mathbb{E} \int_0^t Z_s(m(X_s, \theta_s) + m(X_s, 0))\, ds \leq 2Ct$$

so $(Z_t)_{t \geq 0}$ is a martingale.

PROPOSITION 4.1. *For all $\delta \in (0, A\eta t \sqrt{d}\,]$,*

$$\mathbb{P}\Big( \sup_{s \leq t} |M_s| > \delta \Big) \leq (2d) e^{-\delta^2/(2Adt)}.$$

PROOF.   Fix $\theta \in (\mathbb{R}^d)^*$ with $|\theta| = 1$ and consider the stopping time

$$T = \inf\{t \geq 0 : \langle \theta, M_t \rangle > \delta\}.$$

For $\varepsilon < \eta$, taking $\theta_t = \theta$ for all $t$ above, we know that $(Z_t^{\varepsilon\theta})_{t \geq 0}$ is a martingale. On the set $\{T \leq t\}$ we have $Z_T^{\varepsilon\theta} \geq e^{\delta\varepsilon - At\varepsilon^2/2}$. By optional stopping,

$$\mathbb{E}(Z_{T \wedge t}^{\varepsilon\theta}) = \mathbb{E}(Z_0^{\varepsilon\theta}) = 1.$$

Hence,

$$\mathbb{P}\left(\sup_{s \leq t} \langle \theta, M_s \rangle > \delta\right) = \mathbb{P}(T \leq t) \leq e^{-\delta\varepsilon + At\varepsilon^2/2}.$$

When $\delta \leq At\eta$ we can take $\varepsilon = \delta/At$ to obtain

$$\mathbb{P}\left(\sup_{s \leq t} \langle \theta, M_s \rangle > \delta\right) \leq e^{-\delta^2/2At}.$$

Finally, if $\sup_{s \leq t} |M_s| > \delta$, then $\sup_{s \leq t} \langle \theta, M_s \rangle > \delta/\sqrt{d}$ for one of $\theta = \pm e_1, \ldots, \pm e_d$. □

## 5. Fluid limit for stopped processes.
In this section we develop some general criteria for the convergence of a sequence of Markov chains in $\mathbb{R}^d$ to the solution of a differential equation, paying particular attention to the case where the chain may stop abruptly on leaving a given open set.

5.1. *Fluid limits.* It is possible to give criteria for the convergence of Markov processes in terms of the limiting behavior of their infinitesimal characteristics. This is a powerful technique which has been intensively studied by probabilists. The book of Ethier and Kurtz [8] is a key reference. Further results are given in Chapter 17 of [10] and in Theorems IX.4.21 and IX.4.26 of [9]. A particular case with many applications is where the limiting process is deterministic and is given by a differential equation, sometimes called a fluid limit. The relevant probabilistic literature, though well developed, may not be readily accessible to nonspecialists seeking to apply the results in other fields. One field where fluid limits of Markov processes are beginning to find interesting applications is random combinatorics. Wormald [16] and co-workers have put forward a set of criteria which is specially adapted to this application. The material in this section may be considered as an alternative framework, somewhat more rigid but, we hope, easy to use, developed with the same applications in mind.

Let $(X_t^N)_{t \geq 0}$ be a sequence of pure jump Markov processes in $\mathbb{R}^d$. It may be that $(X_t^N)_{t \geq 0}$ takes values in some discrete subset $I^N$ of $\mathbb{R}^d$ and that its Lévy kernel $K^N(x, dy)$ is given naturally only for $x \in I^N$. So let us suppose that $I^N$ is measurable, that $(X_t^N)_{t \geq 0}$ takes values in $I^N$ and that the Lévy kernel $K^N(x, dy)$ is given for $x \in I^N$. Let $S$ be an open set in $\mathbb{R}^d$ and set $S^N = I^N \cap S$. We shall study, under certain hypotheses, the limiting behavior of $(X_t^N)_{t \geq 0}$ as $N \to \infty$, on compact time intervals, up to the first time the process leaves $S$. In applications, the set $S$ will be chosen as the intersection of two open sets $H$ and $U$. Our sequence of processes may all stop abruptly on leaving some open set $H$, so that $K^N(x, dy) = 0$ for $x \notin H$. If this sort of behavior does not occur, we simply take $H = \mathbb{R}^d$. We choose $U$ so that the conjectured fluid limit path does not leave $U$ in the relevant compact time interval. Subject to this restriction we are free to take $U$ as small as we like to facilitate the checking of convergence and regularity conditions, which are required only on $S$.

The scope of our study is motivated by the particular model which occupies the remainder of this paper: so we are willing to impose a relatively strong, large-deviations-type hypothesis on the Lévy kernels $K^N$, see (6) below, and we are interested to find that strong conclusions may be drawn using rather direct arguments. On the other hand, in certain cases of our model, the fluid limit path grazes the boundary of the set $S$; this calls for a refinement of the usual fluid limit results to determine the limiting distribution of the exit time.

5.2. *Assumptions.* Consider the Laplace transform

$$m^N(x, \theta) = \int_{\mathbb{R}^d} e^{\langle \theta, y \rangle} K^N(x, dy), \qquad x \in S^N, \theta \in (\mathbb{R}^d)^*.$$

We assume that, for some $\eta_0 > 0$,

$$(6) \qquad \sup_N \sup_{x \in S^N} \sup_{|\theta| \leq \eta_0} \frac{m^N(x, N\theta)}{N} < \infty.$$

Set $b^N(x) = m^{N\prime}(x, 0)$, where "$\prime$" denotes the derivative in $\theta$. We assume that, for some Lipschitz vector field $b$ on $S$,

$$(7) \qquad \sup_{x \in S^N} |b^N(x) - b(x)| \to 0.$$

We write $\tilde{b}$ for some Lipschitz vector field on $\mathbb{R}^d$ extending $b$. [Such a $\tilde{b}$ is given, e.g., by $\tilde{b}(x) = \sup\{b(y) - K|x - y| : y \in S\}$ where $K$ is the Lipschitz constant for $b$.] Fix a point $x_0$ in the closure $\bar{S}$ of $S$ and denote by $(x_t)_{t \geq 0}$ the unique solution to $\dot{x}_t = \tilde{b}(x_t)$ starting from $x_0$. We assume finally that, for all $\delta > 0$,

$$(8) \qquad \limsup_{N \to \infty} N^{-1} \log \mathbb{P}(|X_0^N - x_0| > \delta) < 0.$$

While these are not the weakest conditions for the fluid limit, they are readily verified in many examples of interest. In particular, we will be able to verify them for the Markov chains associated with hypergraph collapse in Section 3.

5.3. *Exponential convergence to the fluid limit.* Fix $t_0 > 0$ and set

$$T^N = \inf\{t \geq 0 : X_t^N \notin S\} \wedge t_0.$$

PROPOSITION 5.1.   *Under assumptions* (6)–(8), *we have, for all* $\delta > 0$,

$$\limsup_{N \to \infty} N^{-1} \log \mathbb{P}\left(\sup_{t \leq T^N} |X_t^N - x_t| > \delta\right) < 0. \tag{9}$$

PROOF.   The following argument is widely known but we have not found a convenient reference. Set $b^N(x) = m^{N\prime}(x, 0)$ and define $(M_t^N)_{t \geq 0}$ by

$$X_t^N = X_0^N + M_t^N + \int_0^t b^N(X_s^N)\,ds.$$

Note that $M_t^N$ corresponds to the martingale we identified in Proposition 4.1. Fix $\eta \in (0, \eta_0)$. Assumption (6) implies that there exists $C < \infty$ such that, for all $N$,

$$|m^{N\prime\prime}(x, \theta)| \leq C/N, \qquad x \in S^N, |\theta| \leq N\eta.$$

Compare this estimate with (4). By applying Proposition 4.1 to the stopped process $(X_{t \wedge T^N}^N)_{t \geq 0}$, we find constants $\varepsilon_0 > 0$ and $C_0 < \infty$, depending only on $C, \eta, d$ and $t_0$ such that, for all $N$ and all $\varepsilon \in (0, \varepsilon_0]$,

$$\mathbb{P}\left(\sup_{t \leq T^N} |M_t^N| > \varepsilon\right) \leq C_0 e^{-N\varepsilon^2/C_0}. \tag{10}$$

Given $\delta > 0$, set $\varepsilon = \min\{e^{-Kt_0}\delta/3, \varepsilon_0\}$, where $K$ is the Lipschitz constant of $\tilde{b}$. Let

$$\Omega^N = \left\{|X_0^N - x_0| \leq \varepsilon \text{ and } \sup_{t \leq T^N} |M_t^N| \leq \varepsilon\right\}.$$

Then (8) and (10) together imply that

$$\limsup_{N \to \infty} N^{-1} \log \mathbb{P}(\Omega \backslash \Omega^N) < 0.$$

On the other hand, by (7), there exists $N_0$ such that $|b^N(x) - b(x)| \leq \varepsilon/t_0$ for all $x \in S^N$ and all $N \geq N_0$. We note that

$$X_t^N - x_t = (X_0^N - x_0) + M_t^N + \int_0^t (b^N(X_s^N) - b(X_s^N))\,ds + \int_0^t (\tilde{b}(X_s^N) - \tilde{b}(x_s))\,ds$$

so, for $N \geq N_0$, on $\Omega^N$, for $t \leq T^N$,

$$|X_t^N - x_t| \leq 3\varepsilon + K \int_0^t |X_s^N - x_s|\,ds,$$

which implies, by Gronwall's lemma, that $\sup_{t \leq T^N} |X_t^N - x_t| \leq \delta$.   □

5.4. *Limiting distribution of the exit time.* The remainder of this section is concerned with the question, left open by Proposition 5.1, of determining the limiting distribution of $T^N$. Set

$$\tau = \inf\{t \geq 0 : x_t \notin \bar{S}\} \wedge t_0,$$

$$\mathcal{T} = \{t \in [0, \tau) : x_t \notin S\}.$$

It is straightforward to deduce from (9) that, for all $\delta > 0$,

$$(11) \qquad \limsup_{N \to \infty} N^{-1} \log \mathbb{P}\left(\inf_{t \in \mathcal{T} \cup \{\tau\}} |T^N - t| > \delta\right) < 0.$$

In particular, if $\mathcal{T}$ is empty, then $T^N \to \tau$ in probability and, for all $\delta > 0$,

$$\limsup_{N \to \infty} N^{-1} \log \mathbb{P}\left(\sup_{t \leq t_0} |X_t^N - x_{t \wedge \tau}| > \delta\right) < 0.$$

The reader who wishes only to know the proof of Theorem 2.1 may skip to Section 6 as the remaining results of this section are needed only for the more general case considered in Theorem 2.2.

5.5. *Fluctuations.* We assume here that

$$(12) \qquad \qquad \mathcal{T} \text{ is finite.}$$

In this case the limiting distribution of $T^N$ may be obtained from that of the fluctuations $\gamma_t^N = \sqrt{N}(X_{t \wedge T^N}^N - x_{t \wedge T^N})$. We assume that there exists a limit kernel $K(x, dy)$, defined for $x \in S$, such that $m(x, \theta) < \infty$ for all $x \in S$ and $|\theta| \leq \eta_0$, where

$$m(x, \theta) = \int_{\mathbb{R}^d} e^{\langle \theta, y \rangle} K(x, dy), \qquad x \in S, \theta \in (\mathbb{R}^d)^*.$$

For convergence of the fluctuations we assume

$$(13) \qquad \qquad \gamma_0^N \to \gamma_0 \text{ in distribution,}$$

$$(14) \qquad \sup_{x \in S^N} \sup_{|\theta| \leq \eta_0} \left| \frac{m^N(x, N\theta)}{N} - m(x, \theta) \right| \to 0,$$

$$(15) \qquad \sup_{x \in S^N} \sqrt{N} |b^N(x) - b(x)| \to 0,$$

$$(16) \qquad \qquad a \text{ is Lipschitz and } b \text{ is } C^1 \text{ on } S,$$

where $b^N(x) = m^{N\prime}(x, 0)$ and $a(x) = m''(x, 0)$. Of course (14) will force $b(x) = m'(x, 0)$.

5.6. *Limiting stochastic differential equation.* Consider the process $(\gamma_t)_{t \leq \tau}$ given by the linear stochastic differential equation

$$(17) \qquad d\gamma_t = \sigma(x_t)\, dB_t + \nabla b(x_t)\gamma_t\, dt$$

and starting from $\gamma_0$, where $B$ is a Brownian motion and $\sigma(x)\sigma(x)^* = a(x)$. The distribution of $(\gamma_t)_{t \leq \tau}$ does not depend on the choice of $\sigma$. For convergence of $T^N$ we assume, in addition,

$$(18) \qquad \begin{aligned} &\partial S \text{ is } C^1 \text{ at } x_t \text{ with inward normal } n_t \quad \text{and} \\ &\mathbb{P}(\langle n_t, \gamma_t \rangle = 0) = 0 \qquad \text{for all } t \in \mathcal{T}. \end{aligned}$$

THEOREM 5.1. *Under assumptions* (8), (12)–(16) *and* (18) *we have* $T^N \to T$ *in distribution, where*

$$T = \min\{t \in \mathcal{T} : \langle n_t, \gamma_t \rangle < 0\} \wedge \tau.$$

PROOF. Let $\tau_0 = 0$ and write the positive elements of $\mathcal{T}$ as $\tau_1 < \cdots < \tau_m$. Define, for $k = 0, 1, \ldots, m$,

$$\tilde{\gamma}_k^N = \begin{cases} \gamma_{\tau_k}^N, & \text{if } T^N > \tau_k, \\ \partial, & \text{otherwise}, \end{cases}$$

$$\tilde{\gamma}_k = \begin{cases} \gamma_{\tau_k}, & \text{if } T > \tau_k, \\ \partial, & \text{otherwise}, \end{cases}$$

where $\partial$ is some cemetery state. We will show by induction, for $k = 0, 1, \ldots, m$, that

$$(19) \qquad (\tilde{\gamma}_0^N, \ldots, \tilde{\gamma}_k^N) \to (\tilde{\gamma}_0, \ldots, \tilde{\gamma}_k) \text{ in distribution}.$$

Given (11), this implies that $T^N \to T$ in distribution, as required.

Note that both $(\tilde{\gamma}_k^N)_{0 \leq k \leq m}$ and $(\tilde{\gamma}_k)_{0 \leq k \leq m}$ may be considered as time-dependent Markov processes. Hence, by a conditioning argument, it suffices to deal with the case where $\gamma_0$ is nonrandom. By (18), if $x_0 \in \partial S$, we can assume that $\partial S$ is $C^1$ at $x_0$ and $\langle n_0, \gamma_0 \rangle \neq 0$. Moreover, for the inductive step, it suffices to consider the case where $\tilde{\gamma}_k$ is nonrandom, not $\partial$, and to show that, if $\tilde{\gamma}_k^N \to \tilde{\gamma}_k$ in probability, then $\tilde{\gamma}_{k+1}^N \to \tilde{\gamma}_{k+1}$ in distribution. We lose no generality in considering only the case $k = 0$.

We have assumed that $\gamma_0^N \to \gamma_0$ in distribution. Note that $T = 0$ if and only if $x_0 \in \partial S$ and $\langle n_0, \gamma_0 \rangle < 0$. On the other hand, since $X_0^N = x_0 + \sqrt{N}\gamma_0^N$, we have $\mathbb{P}(T^N = 0) \to 1$ if and only if $x_0 \in \partial S$ and $\langle n_0, \gamma_0 \rangle < 0$. Hence $\tilde{\gamma}_0^N \to \tilde{\gamma}_0$ in distribution; that is, (19) holds for $k = 0$.

In Lemmas 5.4–5.6 below, we will show that, if $x_0 \in S$, or $x_0 \in \partial S$ and $\langle n_0, \gamma_0 \rangle > 0$, then

$$\mathbb{P}(T^N > \varepsilon) \to 1 \qquad \text{for some } \varepsilon > 0,$$

and, in the case $m \geq 1$,

$$\gamma_{\tau_1}^N \to \gamma_{\tau_1} \text{ in distribution,}$$

$$\mathbb{P}(\langle n_{\tau_1}, \gamma_{\tau_1}^N \rangle \geq 0 \text{ and } T^N \leq \tau_1) \to 0,$$

$$\mathbb{P}(\langle n_{\tau_1}, \gamma_{\tau_1}^N \rangle < 0 \text{ and } T^N > \tau_1) \to 0.$$

It follows that $\tilde{\gamma}_1^N \to \tilde{\gamma}_1$ in distribution, so (19) holds for $k = 1$. This establishes the induction and completes the proof. $\square$

We remark that the same proof applies when the Lévy kernels $K^N$ have a measurable dependence on the time parameter $t$, subject to obvious modifications and to each hypothesis holding uniformly in $t \leq t_0$.

For the remainder of this section, the assumptions of Theorem 5.1 are in force and $\gamma_0$ is nonrandom.

LEMMA 5.2. *For all $\varepsilon > 0$ there exists $\lambda < \infty$ such that, for all $N$,*

$$\mathbb{P}\left(\sup_{t \leq t_0} |\gamma_t^N| \geq \lambda\right) < \varepsilon.$$

PROOF. Given $\varepsilon > 0$, choose $\lambda < \infty$ and $N_0$ such that, for $\lambda' = e^{-Kt_0}\lambda/3$ and $N \geq N_0$,

$$\sqrt{N}|b^N(x) - b(x)| \leq \lambda'/t_0, \qquad x \in S^N,$$

and, with probability exceeding $1 - \varepsilon$,

$$|\gamma_0^N| \leq \lambda',$$

$$\sqrt{N} \sup_{t \leq T^N} |M_t^N| \leq \lambda'.$$

This is possible by (10) and (15). These three inequalities imply

$$|\gamma_t^N| \leq 3\lambda' + K \int_0^t |\gamma_s^N| \, ds, \qquad t \leq T^N,$$

so, by Gronwall's lemma

$$\sup_{t \leq T^N} |\gamma_t^N| \leq \lambda. \qquad \square$$

LEMMA 5.3. *For all $\varepsilon > 0$ there exists $\lambda < \infty$ such that, for all $\delta > 0$, there exists $N_\delta < \infty$ such that, for all $N \geq N_\delta$ and all $t \leq t_0$,*

$$\mathbb{P}\left(\sup_{s \leq t_0, t \leq s \leq t+\delta} |\gamma_s^N - \gamma_t^N| > \lambda\sqrt{\delta}\right) < \varepsilon.$$

PROOF. Consider first the case $t = 0$. Given $\varepsilon > 0$, choose $\lambda < \infty$ such that, for all $\delta > 0$, there exists $N_\delta$ such that, for $\lambda' = e^{Kt_0}\lambda/3$ and $N \geq N_\delta$,

$$\sqrt{N}|b^N(x) - b(x)| \leq \lambda'/\sqrt{t_0}, \qquad x \in S^N,$$

and, with probability exceeding $1 - \varepsilon$,

$$|\gamma_0^N| \leq \lambda'/K\sqrt{t_0},$$

$$\sqrt{N} \sup_{t \leq T^N \wedge \delta} |M_t^N| \leq \lambda'\sqrt{\delta}.$$

This is possible by (10) and (15). These three inequalities imply

$$|\gamma_t^N - \gamma_0^N| \leq 3\lambda'\sqrt{\delta} + K\int_0^t |\gamma_s^N - \gamma_0^N|\, ds, \qquad t \leq T^N \wedge \delta,$$

so by Gronwall's lemma,

$$\sup_{t \leq T^N \wedge \delta} |\gamma_t^N - \gamma_0^N| \leq \lambda\sqrt{\delta}.$$

The case $t > 0$ follows by the same sort of argument, using Lemma 5.2 to get the necessary tightness of $\gamma_t^N$. $\square$

LEMMA 5.4. *Suppose either $x_0 \in S$, or $x_0 \in \partial S$ and $\langle n_0, \gamma_0 \rangle > 0$. Then $\mathbb{P}(T^N > \varepsilon) \to 1$ as $N \to \infty$ for some $\varepsilon > 0$.*

PROOF. The case $x_0 \in S$ follows from (11). Suppose then that $x_0 \in \partial S$ and $\langle n_0, \gamma_0 \rangle > 0$. Then, since $\partial S$ is $C^1$ at $x_0$, for all $\varepsilon > 0$, there exists $\delta(\varepsilon) > 0$ such that, for all $x \in \bar{S}$ with $|x - x_0| \leq \delta(\varepsilon)$, and all $v \in \mathbb{R}^d$,

$$(20) \qquad |v| \leq \delta(\varepsilon) \quad \text{and} \quad \langle n_0, v \rangle \geq \varepsilon|v| \quad \Longrightarrow \quad x + v \in S.$$

Since $\langle n_0, \gamma_0 \rangle > 0$, by Lemma 5.3, given $\varepsilon > 0$ there exist $\varepsilon_1 > 0$ and $N_0$ such that, for all $N \geq N_0$ and $t \leq T^N \wedge \varepsilon_1$,

$$\langle n_0, \gamma_t^N \rangle > \varepsilon_1|\gamma_t^N|, \qquad |\gamma_t^N| < 1/\varepsilon_1,$$

with probability exceeding $1 - \varepsilon$. Choose $\varepsilon_2 \in (0, \varepsilon_1)$ so that $|x_t - x_0| \leq \delta(\varepsilon_1)$ and $x_t \in \bar{S}$ whenever $t \leq \varepsilon_2$. Set $N_1 = \max\{N_0, (\varepsilon_1\delta(\varepsilon_1))^{-2}\}$, then, for $N \geq N_1$ and $t \leq T^N \wedge \varepsilon_2$,

$$(21) \qquad x_t \in \bar{S}, \qquad |x_t - x_0| \leq \delta(\varepsilon_1), \qquad N^{-1/2}|\gamma_t^N| \leq \delta(\varepsilon_1),$$
$$\langle n_0, \gamma_t^N \rangle > \varepsilon_1|\gamma_t^N|,$$

with probability exceeding $1 - \varepsilon$. By (20), (21) implies $X_t^N = x_t + N^{-1/2}\gamma_t^N \in S$. Hence $\mathbb{P}(T^N \leq \varepsilon_2) < \varepsilon$ for all $N \geq N_1$. $\square$

For the rest of this section we assume that $m \geq 1$. (The next result holds with $\tau_1$ replaced by $\tau$ when $m = 0$, by the same argument, but we do not need this.)

LEMMA 5.5. *Suppose either $x_0 \in S$, or $x_0 \in \partial S$ and $\langle n_0, \gamma_0 \rangle > 0$. Then $\gamma_{\tau_1}^N \to \gamma_{\tau_1}$ in distribution as $N \to \infty$.*

PROOF. By Lemma 5.3, given $\delta > 0$, we can find $t < \tau_1$ such that, for all $N$,

$$\mathbb{P}(|\gamma_t^N - \gamma_{\tau_1}^N| > \delta) < \delta, \qquad \mathbb{P}(|\gamma_t - \gamma_{\tau_1}| > \delta) < \delta.$$

Hence it suffices to show $\gamma_t^N \to \gamma_t$ in distribution for all $t < \tau_1$.

Define $(\psi_t)_{t \le \tau}$ in $\mathbb{R}^d \otimes (\mathbb{R}^d)^*$ by

$$\dot\psi_t = \nabla b(x_t)\psi_t, \qquad \psi_0 = id.$$

Fix $\theta \in (\mathbb{R}^d)^*$ and set $\theta_t = (\psi_t^*)^{-1}\theta$. Then

$$d\langle \theta_t, \gamma_t \rangle = \langle \theta_t, \sigma(x_t)\,dB_t \rangle, \qquad t \le \tau,$$

so

$$\langle \theta_t, \gamma_t \rangle \sim N\left(\langle \theta, \gamma_0 \rangle, \int_0^t \langle \theta_s, a(x_s)\theta_s \rangle\,ds\right), \qquad t \le \tau.$$

On the other hand, for $(M_t^N)_{t \ge 0}$ as in the proof of Proposition 5.1,

$$d\langle \theta_t, \gamma_t^N \rangle = \sqrt{N}\langle \theta_t, dM_t^N \rangle + R_t^{N,\theta}\,dt, \qquad t \le T^N,$$

where

$$R_t^{N,\theta} = \sqrt{N}\langle \theta_t, b^N(X_t^N) - b(x_t) - \nabla b(x_t)(X_t^N - x_t) \rangle.$$

By (15),

$$\sup_{t \le T^N} \sqrt{N}|b^N(X_t^N) - b(X_t^N)| \to 0.$$

By (16), given $\varepsilon > 0$, there exists $\delta > 0$ such that, for all $t \in [\varepsilon, \tau_1 - \varepsilon]$, for $|x - x_t| \le \delta$,

$$|b(x) - b(x_t) - \nabla b(x_t)(x - x_t)| \le \varepsilon|x - x_t|.$$

Hence $|X_t^N - x_t| \le \delta$ and $\varepsilon \le t \le \tau_1 - \varepsilon$ imply

$$\sqrt{N}|b(X_t^N) - b(x_t) - \nabla b(x_t)(X_t^N - x_t)| \le \varepsilon|\gamma_t^N|.$$

Combining this with Lemma 5.2, we deduce that

$$\int_0^{\tau_1} |R_t^{N,\theta}|\,dt \to 0 \text{ in probability.}$$

Hence it suffices to show, for all $\theta \in (\mathbb{R}^d)^*$ and all $t < \tau_1$,

$$\sqrt{N}\int_0^t \langle \theta_s, dM_s^N \rangle \to N\left(0, \int_0^t \langle \theta_s, a(x_s)\theta_s \rangle\,ds\right) \text{ in distribution.}$$

Indeed, it suffices to show, for all $\theta \in (\mathbb{R}^d)^*$ and $t < \tau_1$, that $\mathbb{E}(E_t^{N,\theta}) \to 1$ as $N \to \infty$, where

$$E_t^{N,\theta} = \exp\left\{ i\sqrt{N} \int_0^t \langle \theta_s, dM_s^N \rangle + \tfrac{1}{2} \int_0^t \langle \theta_s, a(x_s)\theta_s \rangle \, ds \right\}.$$

Set $\tilde{m}^N(x, \theta) = m^N(x, i\theta)$, $\tilde{m}(x, \theta) = m(x, i\theta)$ and

$$\tilde{\phi}^N(x, \theta) = \int_{\mathbb{R}^d} (e^{i\langle \theta, y \rangle} - 1 - i\langle \theta, y \rangle) K^N(x, dy).$$

By (14), for all $\eta < \eta_0$, we have

$$\sup_{x \in S^N} \sup_{|\theta| \le \eta} |N\tilde{m}^{N\prime\prime}(x, N\theta) - \tilde{m}''(x, \theta)| \to 0.$$

Note that

$$\tilde{\phi}^N(x, \sqrt{N}\theta) + \tfrac{1}{2}\langle \theta, a(x)\theta \rangle$$
$$= \int_0^1 (N\tilde{m}^{N\prime\prime}(x, \sqrt{N}r\theta) - \tilde{m}''(x, 0))(\theta, \theta)(1 - r) \, dr$$

so, for all $\rho < \infty$,

(22)     $$\sup_{x \in S^N} \sup_{|\theta| \le \rho} |\tilde{\phi}^N(x, \sqrt{N}\theta) + \tfrac{1}{2}\langle \theta, a(x)\theta \rangle| \to 0.$$

Write $E_t^{N,\theta} = E_t^N = Z_t^N A_t^N B_t^N$, where

$$Z_t^N = \exp\left\{ i\sqrt{N} \int_0^t \langle \theta_s, dM_s^N \rangle - \int_0^t \tilde{\phi}^N(X_s^N, \sqrt{N}\theta_s) \, ds \right\},$$

$$A_t^N = \exp\left\{ \int_0^t (\tilde{\phi}^N(X_s^N, \sqrt{N}\theta_s) + \tfrac{1}{2}\langle \theta_s, a(X_s^N)\theta_s \rangle) \, ds \right\},$$

$$B_t^N = \exp\left\{ \int_0^t \tfrac{1}{2}\langle \theta_s, (a(x_s) - a(X_s^N))\theta_s \rangle \, ds \right\}.$$

Now $(Z_{t \wedge T^N}^N)_{t \le \tau}$ is a martingale, as in (5), so $\mathbb{E}(Z_{t \wedge T^N}^N) = 1$ for all $N$. Fix $t \le \tau$. By (22), $Z_{t \wedge T^N}^N$ is bounded, uniformly in $N$, and $A_{t \wedge T^N}^N \to 1$ uniformly as $N \to \infty$. Moreover, by (16), $B_{t \wedge T^N}^N$ is bounded uniformly in $N$ and converges to 1 in probability, using (9). Hence

$$\mathbb{E}(Z_{t \wedge T^N}^N A_{t \wedge T^N}^N B_{t \wedge T^N}^N) \to 1$$

as $N \to \infty$. By Lemma 5.4 and (11), $\mathbb{P}(T^N > t) \to 1$ for all $t < \tau_1$. It follows that $\mathbb{E}(E_t^N) \to 1$ for all $t < \tau_1$ as required.   $\square$

LEMMA 5.6.   *Suppose either $x_0 \in S$, or $x_0 \in \partial S$ and $\langle n_0, \gamma_0 \rangle > 0$. Then, as $N \to \infty$,*

$$\mathbb{P}(\langle n_{\tau_1}, \gamma_{\tau_1}^N \rangle \ge 0 \ and \ T^N \le \tau_1) \to 0,$$

$$\mathbb{P}(\langle n_{\tau_1}, \gamma_{\tau_1}^N \rangle < 0 \ and \ T^N > \tau_1) \to 0.$$

Proof. By Lemma 5.5, given $\varepsilon > 0$, there exists $\varepsilon_1 > 0$ and $N_0$ such that, for all $N \geq N_0$,

$$|\langle n_{\tau_1}, \gamma_{\tau_1}^N \rangle| > \varepsilon_1 |\gamma_{\tau_1}^N|, \qquad \varepsilon_1 < |\gamma_{\tau_1}^N| < 1/\varepsilon_1,$$

with probability exceeding $1 - \varepsilon$. Then by Lemma 5.3, there exists $\varepsilon_2 > 0$ and $N_1 \geq N_0$ such that, for all $N \geq N_1$, with probability exceeding $1 - \varepsilon$, *either*

$$(23) \qquad \langle n_{\tau_1}, \gamma_t^N \rangle > \varepsilon_2 |\gamma_t^N|, \qquad |\gamma_t^N| < 1/\varepsilon_2 \qquad \text{for all } t \in [\tau_1 - \varepsilon_2, \tau_1]$$

or

$$(24) \qquad \langle n_{\tau_1}, \gamma_{\tau_1}^N \rangle < -\varepsilon_2 |\gamma_{\tau_1}^N|, \qquad |\gamma_{\tau_1}^N| < 1/\varepsilon_2.$$

Since $\partial S$ is $C^1$ at $x_{\tau_1}$, there exists $\delta > 0$ such that

$$\text{if } x \in \bar{S} \text{ and } v \in \mathbb{R}^d \text{ with } |x - x_{\tau_1}| \leq \delta, |v| \leq \delta \quad \text{and}$$

$$\langle n_{\tau_1}, v \rangle < \varepsilon_2 |v| \qquad \text{then } x + v \in S$$

and

$$\text{if } v \in \mathbb{R}^d \text{ with } |v| \leq \delta \quad \text{and} \quad \langle n_{\tau_1}, v \rangle < -\varepsilon_2 |v| \qquad \text{then } x_{\tau_1} + v \notin S.$$

Choose $\varepsilon_3 \in (0, \varepsilon_2]$ such that $|x_t - x_{\tau_1}| \leq \delta$ and $x_t \in \bar{S}$ whenever $t \in [\tau_1 - \varepsilon_3, \tau_1]$. Set $N_2 = \max\{N_1, (\varepsilon_2\delta)^{-2}\}$. Then, for $N \geq N_2$, since $X_t^N = x_t + N^{-1/2}\gamma_t^N$ on $\{T^N \geq t\}$, (23) implies $X_t^N \in S$ for all $t \in [\tau_1 - \varepsilon_2, \tau_1]$ or $T^N < \tau_1 - \varepsilon_2$, and (24) implies $X_{\tau_1}^N \notin S$ or $T^N < \tau_1$. We know by Lemma 5.4 and (11) that $\mathbb{P}(T^N < \tau_1 - \varepsilon_2) \to 0$ as $N \to \infty$. Hence, with high probability, as $N \to \infty$, $\langle n_{\tau_1}, \gamma_{\tau_1}^N \rangle \geq 0$ implies (23) and then $T^N > \tau_1$, and $\langle n_{\tau_1}, \gamma_{\tau_1}^N \rangle < 0$ implies (24) and then $T^N \leq \tau_1$. $\square$

**6. Fluid limit of collapsing hypergraphs.** We now apply the general theory from the preceding sections to prove our main results Theorems 2.1 and 2.2.

6.1. *Lévy kernel for collapse of random hypergraphs.* In Section 3 we introduced a Markov process $(\Lambda_n)_{n \geq 0}$ of collapsing hypergraphs, starting from $\Lambda_0 \sim \text{Poisson}(\beta)$ and stopping when $n = |V^*|$, the number of identifiable vertices in $\Lambda_0$. The process $(Y_n, Z_n)_{n \geq 0}$ of patches and debris in $\Lambda_n$ was found itself to be Markov. We now view this process as a function of the initial number of vertices $N$ and obtain a fluid limit result when $N \to \infty$.

It will be convenient to embed our process in continuous time, by removing vertices according to a Poisson process $(\nu_t)_{t \geq 0}$ of rate $N$ which stops when $\nu_t = |V^*|$. Set

$$X_t^N = N^{-1}(\nu_t, Y_{\nu_t}, Z_{\nu_t})$$

and note that $X^N$ takes values in

$$
(25) \quad \begin{aligned}
I^N = \{x \in \mathbb{R}^3 : Nx^1 \in \{0, 1, \dots, N-1\}, Nx^2, Nx^3 \in \mathbb{Z}^+\} \\
\cup \{(1, 0, x^3) : Nx^3 \in \mathbb{Z}^+\}.
\end{aligned}
$$

The Lévy kernel $K^N(x, dy)$ for $(X_t^N)_{t \geq 0}$ is naturally defined for $x \in I^N$. If $x^2 = 0$, then $K^N(x, dy) = 0$. If $x^2 > 0$, then $N^{-1}K^N(x, \cdot)$ is a probability measure; by Lemma 3.1, it is the law of the random variable $J^N/N$, where

$$
J^N = (1, -1 - W^N + U^N, 1 + W^N),
$$

$$
W^N \sim B(Nx^2 - 1, 1/(N - Nx^1)), \qquad U^N \sim P((N - Nx^1 - 1)\lambda_2(N, Nx^1))
$$

with $W^N$ and $U^N$ independent.

Recall that $R$ denotes the radius of convergence of the power series $\beta(t)$, given by (1). We assume, until further notice, that $R > 0$ and fix $t_0 \in (0, R \wedge 1)$ and $\rho \in (t_0, R \wedge 1)$.

LEMMA 6.1.  *There is a constant $C < \infty$ such that*

$$
|N\lambda_2(N, n) - \beta''(n/N)| \leq C(\log N)^2/N
$$

*for all $N \in \mathbb{N}$ and $n \in \{0, 1, \dots, [N\rho]\}$.*

PROOF.   Recall that

$$
\lambda_2(N, n) = N \sum_{i=0}^{n} (i+1)(i+2)\beta_{i+2} \frac{n(n-1)\cdots(n-i+1)}{N(N-1)\cdots(N-i-1)}.
$$

Set $M = A \log N$ where $A = (\log(R/\rho))^{-1} < \infty$. Then, for $n \leq [N\rho]$,

$$
\begin{aligned}
&|N\lambda_2(N, n) - \beta''(n/N)| \\
&\qquad \leq \sum_{i=1}^{M \wedge n} (i+1)(i+2)\beta_{i+2}\delta_i(N, n) + 2 \sum_{i=M+1}^{\infty} (i+1)(i+2)\beta_{i+2}\rho_i(N, n),
\end{aligned}
$$

where

$$
\delta_i(N, n) = \left| \frac{N^2}{(N-i)(N-i-1)} \frac{n}{N} \left( \frac{n-1}{N-1} \right) \cdots \left( \frac{n-i+1}{N-i+1} \right) - \left( \frac{n}{N} \right)^i \right|
$$

and

$$
\rho_i(N, n) = \frac{N^2}{(N-i)(N-i-1)} \left( \frac{n}{N} \right)^i \leq C\rho^i.
$$

Note that, for $j = 0, \dots, i-1$ and $i \leq M \wedge n$,

$$
\left| \frac{n-j}{N-j} - \frac{n}{N} \right| \leq A \log N/N
$$

so, making use of the inequality $|\prod a_j - \prod b_j| \leq \sum |a_j - b_j|$ for $0 \leq a_j, b_j \leq 1$, we obtain

$$\delta_i(N, n) \leq C(\log N)^2 \rho^i / N.$$

Hence

$$|N\lambda_2(N, n) - \beta''(n/N)| \leq C(\log N)^2 \beta''(\rho)/N + C(\rho/R)^M$$

and $(\rho/R)^M = 1/N$. $\quad\square$

6.2. *Fluid limit.* The main result of this section is to obtain the limiting behavior of $(X_t^N)_{t \geq 0}$ as $N \to \infty$, which we deduce from Proposition 5.1 and Theorem 5.1. We present first the calculations by which the limit was discovered.

Note that, as $N \to \infty$, for $x^1 < R \wedge 1$, we have $W^N \to W$ and $U^N \to U$ in distribution, where

$$W \sim P(x^2/(1-x^1)), \qquad U \sim P((1-x^1)\beta''(x^1)).$$

Set $J = (1, -1 - W + U, 1 + W)$. Note also that $X_0^N \to x_0 = (0, \beta_1, \beta_0)$ and $\sqrt{N}(X_0^N - x_0) \to \gamma_0$ in distribution, where $\gamma_0^1 = 0, \gamma_0^2 \sim N(0, \beta_1), \gamma_0^3 \sim N(0, \beta_0)$, with $\gamma_0^2$ and $\gamma_0^3$ independent. Thus, subject to certain technical conditions, to be checked later, at least up to the first time that $X_t^{N,1} \geq R \wedge 1$ or $X_t^{N,2} = 0$, the limit path is given by $\dot{x}_t = b(x_t)$, starting from $x_0$, where

$$b(x) = \mathbb{E}(J) = \left(1, -1 - \frac{x^2}{1-x^1} + (1-x^1)\beta''(x^1), \frac{x^2}{1-x^1}\right).$$

Fix $\rho' \in (0, \infty)$ and set

(26) $$S = \{(x^1, x^2, x^3) : |x^1| < \rho, x^2 \in (0, \rho'), x^3 \in \mathbb{R}\}.$$

Then $b$ is Lipschitz on $S$ and, for $\rho'$ sufficiently large, the maximal solution on $[0, t_0]$ to $\dot{x}_t = b(x_t)$ in $\bar{S}$ starting from $x_0$ is given by $(x_t)_{t \leq \tau}$, where

$$x_t = (t, (1-t)(\beta'(t) + \log(1-t)), \beta(t) - (1-t)\log(1-t))$$

and

$$\tau = z^* \wedge t_0.$$

6.3. *Limiting fluctuations.* Set $a(x) = \mathbb{E}(J \otimes J)$. A convenient choice of $\sigma$ such that $\sigma\sigma^* = a$ is $\sigma = (V_1, V_2, V_3)$, where

$$V_1(x) = \sqrt{\frac{x^2}{1-x^1}} \begin{pmatrix} 0 \\ 1 \\ -1 \end{pmatrix},$$

$$V_2(x) = \sqrt{(1-x^1)\beta''(x^1)} \begin{pmatrix} 0 \\ 1 \\ 0 \end{pmatrix},$$

$$V_3(x) = b(x).$$

Note that $a$ is Lipschitz and $b$ is $C^1$ on $S$. The limiting fluctuations are given by

$$d\gamma_t = \sum_i V_i(x_t)\, dB_t^i + \nabla b(x_t)\gamma_t\, dt, \qquad t \le \tau,$$

starting from $\gamma_0$, where $B$ is a Brownian motion in $\mathbb{R}^3$ independent of $\gamma_0$. Note that

$$\mathcal{T} = \{t \in [0,\tau) : x_t \notin S\} = \zeta \cap [0, t_0).$$

In cases where $\mathcal{T}$ is nonempty, the limiting behavior of $(X_t^N)_{t \le t_0}$ depends on the signs of the component of the fluctuations normal to the boundary, that is, on $(\gamma_t^2 : t \in \mathcal{T})$.

Note that $\theta_t = b(x_t)B_t^3$ satisfies

$$d\theta_t = V_3(x_t)\, dB_t^3 + \nabla b(x_t)\theta_t\, dt.$$

This is the part of the fluctuations which reflects our Poissonization of the time-scale. Since $b^2(x_t) = 0$ for all $t \in \mathcal{T}$, it does not affect $(\gamma_t^2 : t \in \mathcal{T})$. So consider $\gamma_t^* = \gamma_t - \theta_t$. Then

$$d\gamma_t^* = V_1(x_t)\, dB_t^1 + V_2(x_t)\, dB_t^2 + \nabla b(x_t)\gamma_t^*\, dt.$$

Note that $V_1^1(x) = V_2^1(x) = 0$ and $\nabla b^1(x) = 0$, so $(\gamma_t^*)^1 = (\gamma_0^*)^1 = 0$ for all $t$. Also $\partial b^2/\partial x^2 = -1/(1 - x^1)$ and $\partial b^2/\partial x^3 = 0$. Also $x_t^1 = t$ and $x_t^2/(1 - x_t^1) = \beta'(t) + \log(1-t)$. The sign of $(\gamma_t^*)^2$ is the same as that of $\alpha_t = (\gamma_t^*)^2/(1-t)$. We have

$$d\alpha_t = d\gamma_t^{*2}/(1-t) + \gamma_t^{*2}/(1-t)^2\, dt$$
$$= (V_1^2(x_t)\, dB_t^1 + V_2^2(x_t)\, dB_t^2)/(1-t)$$

so we can write $\alpha_t = W(\sigma_t^2)$, where $W$ is a Brownian motion and

$$\sigma_t^2 = \beta_1 + \int_0^t \frac{\beta'(s) + \log(1-s) + (1-s)\beta''(s)}{1-s}\, ds$$
$$= \frac{\beta'(t) + \log(1-t) + t}{1-t}.$$

We have shown that $(\mathrm{sgn}(\gamma_t^2) : t \in \mathcal{T})$ has the same distribution as $(\mathrm{sgn}(W_{t/(1-t)}) : t \in \mathcal{T})$. In particular, $\mathbb{P}(\gamma_t^2 = 0) = 0$ for all $t \in \mathcal{T}$.

Recall that $Z$ is defined by

$$Z = \min\{z \in \zeta : W(z/1 - z) < 0\} \wedge z^*.$$

Set

$$T^N = \inf\{t \ge 0 : X_t^{N,2} = 0\}$$

and put $Z(t_0) = Z \wedge t_0$, $T^N(t_0) = T^N \wedge t_0$.

THEOREM 6.1. *For all $\delta > 0$ we have*

$$\limsup_{N \to \infty} N^{-1} \log \mathbb{P}\left(\sup_{t \le T^N(t_0)} |X_t^N - x_t| > \delta\right) < 0.$$

*Moreover, $T^N(t_0) \to Z(t_0)$ in distribution as $N \to \infty$.*

PROOF. We defined $I^N$, the state-space of $(X_t^N)_{t \ge 0}$, in (25), and $S$ in (26). Set $S^N = I^N \cap S$. For $x \in S^N$ we have

$$m^N(x, \theta) = \int_{\mathbb{R}^3} e^{\langle \theta, y \rangle} K^N(x, dy) = N \mathbb{E}(e^{\langle \theta, J^N \rangle / N})$$

so

$$\frac{m^N(x, \theta)}{N} = \exp\left\{\theta_1 - \theta_2 + \theta_3 + B\left(Nx^2 - 1, \frac{1}{N - Nx^1}, \frac{\theta_3 - \theta_2}{N}\right)\right.$$
$$\left. + P\left((N - Nx^1 - 1)\lambda_2(N, Nx^1), \frac{\theta_2}{N}\right)\right\},$$

where, for $\theta \in \mathbb{R}$, we write $B(N, p, \theta) = N \log(1 - p + pe^\theta)$ and $P(\lambda, \theta) = \lambda(e^\theta - 1)$. So, by Lemma 6.1,

$$\sup_{x \in S^N} \sup_{|\theta| \le \eta_0} \left|\frac{m^N(x, N\theta)}{N} - m(x, \theta)\right| \to 0$$

as $N \to \infty$, for all $\eta_0 > 0$, where

$$m(x, \theta) = \mathbb{E}(e^{\langle \theta, J \rangle})$$
$$= \exp\left\{\theta_1 - \theta_2 + \theta_3 + P\left(\frac{x^2}{1 - x^1}, \theta_3 - \theta_2\right) + P((1 - x^1)\beta''(x^1), \theta_2)\right\}.$$

Set

$$b^N(x) = \int_{\mathbb{R}^3} y K^N(x, dy) = \mathbb{E}(J^N),$$

then, by Lemma 7.1,

$$\sup_{x \in S^N} \sqrt{N} |b^N(x) - b(x)| \to 0.$$

Recall that

$$X_0^{N,1} = 0, \qquad NX_0^{N,2} \sim P(N\beta_1), \qquad NX_0^{N,3} \sim P(N\beta_0),$$

and $x_0 = (0, \beta_1, \beta_0)$. By standard exponential estimates, for all $\delta > 0$

$$\limsup_{N \to \infty} N^{-1} \log \mathbb{P}(|X_0^N - x_0| > \delta) < 0.$$

We have now checked the validity of (7), (8), (12)–(16) and (18) in this context, so Proposition 5.1 and Theorem 5.1 apply to give the desired conclusions. □

REMARK 6.1. If $z^* < 1$, then $z^* < R \wedge 1$, so by choosing $t_0 \in (z^*, R \wedge 1)$ we get $Z(t_0) = Z$ and, as $N \to \infty$, with high probability $T^N(t_0) = T^N$. Hence, when $z^* < 1$, Theorem 6.1 holds with $Z$ and $T^N$ replacing $Z(t_0)$ and $T^N(t_0)$. In particular, Theorem 2.1 follows.

6.4. *Proof of Theorem 2.2.* Recall that

$$X_{T^N}^N = \left( \frac{|V^{N*}|}{N}, 0, \frac{|\Lambda^{N*}|}{N} \right).$$

Let $z \in \zeta \cup \{z^*\}$. If $z < 1$, then $z < R \wedge 1$ so, by choosing $t_0 \hat{\in} (z, R \wedge 1)$ in Theorem 5.1, we obtain

$$(27) \qquad \mathbb{P}(|X_{T^N}^N - x_z| \leq \delta) \to \mathbb{P}(Z = z)$$

for all sufficiently small $\delta > 0$.

It remains to deal with the case $z = z^* = 1$. Note that $|V^{N*}| \leq N$ and $|\Lambda^{N*}| \leq |\Lambda^N|$. Now $|\Lambda^N| \sim P(N\beta(1))$ so $|\Lambda^N|/N \to \beta(1)$ in probability as $N \to \infty$. It therefore suffices to show, for all $\delta > 0$ and $\alpha < \beta(1) - \delta$,

$$\liminf_{N \to \infty} \mathbb{P}\left( \frac{|V^{N*}|}{N} \geq 1 - \delta \text{ and } \frac{|\Lambda^{N*}|}{N} \geq \alpha \right) \geq \mathbb{P}(Z = 1).$$

When combined with (27) this completes the proof as we have exhausted the possible values of $Z$.

We consider first the case $R \geq 1$. We can find $t_0 \in (1 - \delta/2, 1)$ such that $\beta(t_0) > \alpha + \delta/2$. Note that $|X_{t_0}^N - x_{t_0}| \leq \delta/2$ implies

$$|V^{N*}|/N \geq X_{t_0}^{N,1} \geq t_0 - \delta/2 > 1 - \alpha,$$

$$|\Lambda^{N*}|/N \geq X_{t_0}^{N,3} \geq \beta(t_0) - (1 - t_0)\log(1 - t_0) - \delta/2 > \alpha.$$

By Theorem 6.1

$$\liminf_{N \to \infty} \mathbb{P}\left( \sup_{t \leq t_0} |X_t^N - x_t| \leq \delta/2 \right) \geq \mathbb{P}(Z > t_0) \geq \mathbb{P}(Z = 1)$$

so we are done.

Consider next the case $R = 0$. Fix $M \in \mathbb{N}$ and set $\tilde{\beta}_j = \beta_j$ if $j \leq M$ and $\tilde{\beta}_j = 0$ otherwise. Then, with obvious notation, we can choose $M$ so that $\tilde{z}_0 > 1 - \delta/2$, $\tilde{\zeta} = \varnothing$ and $\tilde{\beta}(z_0) > \alpha + \delta/2$. Hence

$$\mathbb{P}\left( \frac{|\tilde{V}^{N*}|}{N} \geq 1 - \delta \text{ and } \frac{|\tilde{\Lambda}^{N*}|}{N} \geq \alpha \right) \to 1.$$

We can couple $\Lambda$ and $\tilde{\Lambda}$ so that $\tilde{\Lambda}(A) = \Lambda(A)\mathbb{1}_{|A| \leq M}$. Then $\tilde{V}^{N*} \subseteq V^{N*}$ and $\tilde{\Lambda}^{N*} \leq \Lambda^{N*}$, so this is enough.

There remains the case $R \in (0, 1)$. In this case $\zeta$ is finite. We have assumed that $R \notin \zeta$. So we can find $\rho \in (\sup \zeta, R)$ and $M \in \mathbb{N}$ such that, with obvious notation,

$$\tilde{z}_0 > 1 - \delta/2, \qquad \tilde{\zeta} = \zeta, \qquad \tilde{\beta}(\tilde{z}_0) > \alpha + \delta/2,$$

where $\tilde{\beta}(t), t \in [0, 1)$, is defined by

$$\tilde{\beta}(0) = \beta_0, \qquad \tilde{\beta}'(0) = \beta_1, \qquad \tilde{\beta}''(t) = \begin{cases} \beta''(t), & t < \rho, \\ \sum_{j=2}^{M} j(j-1)\beta_j t^{j-2}, & t \geq \rho. \end{cases}$$

Consider the collapsing hypergraph $(\tilde{\Lambda}_n^N)_{n \geq 0}$ which evolves as $(\Lambda_n^N)_{n \geq 0}$ up to $n = \nu(\rho)$, at which time all hyperedges having at least two vertices and originally having more than $M$ vertices are removed, so that $\tilde{\Lambda}_{\nu(\rho)}^N \leq \Lambda_{\nu(\rho)}^N$. After $\nu(\rho), (\tilde{\Lambda}_n^N)_{n \geq 0}$ evolves by selection of patches as before. Denote by $\tilde{V}^{N*}$ the set of identifiable vertices in $\tilde{\Lambda}_{\nu(\rho)}^N$ and by $\tilde{\Lambda}^{N*}$ the corresponding identifiable hypergraph. Then

$$\tilde{V}^{N*} \subseteq V^{N*} \quad \text{and} \quad \tilde{\Lambda}^{N*} \leq \Lambda^{N*}.$$

A modification of Theorem 5.1 shows that

$$\mathbb{P}(|\tilde{X}_1^N - \tilde{x}_{\tilde{z}_0}| \leq \delta/2) \to \mathbb{P}(\tilde{Z} = \tilde{z}_0) = \mathbb{P}(Z = 1)$$

with $\tilde{X}_1^N = (|\tilde{V}^{N*}|/N, 0, |\tilde{\Lambda}^{N*}|/N)$ and with

$$\tilde{x}_t = (t, (1-t)(\tilde{\beta}'(t) + \log(1-t)), \tilde{\beta}(t) - (1-t)\log(1-t)).$$

All that changes in the proof is that, for $t \geq \rho$ the Lévy kernel is modified by replacing $\lambda_2$ by $\tilde{\lambda}_2$ given by

$$\tilde{\lambda}_2(N, n) = N \sum_{i=0}^{n \wedge (M-2)} (i+1)(i+2)\beta_{i+2} \frac{n(n-1)\cdots(n-i+1)}{N(N-1)\cdots(N-i+1)}.$$

The argument of Lemma 6.1 shows that for all $\rho' < 1$ there is a constant $C < \infty$ such that

$$|N\tilde{\lambda}_2(N, n) - \tilde{\beta}''(n/N)| \leq C/N$$

for all $N \in \mathbb{N}$ and $n = \{0, 1, \ldots, [N\rho']\}$. Everything else is the same. Now $|\tilde{X}_1^N - \tilde{x}_{\tilde{z}_0}| \leq \delta/2$ implies

$$|V^{N*}|/N \geq |\tilde{V}^{N*}|/N = \tilde{X}_1^{N,1} \geq \tilde{z}_0 - \delta/2 \geq 1 - \delta,$$

$$|\Lambda^{N*}|/N \geq |\tilde{\Lambda}^{N*}|/N = \tilde{X}_1^{N,3} \geq \tilde{\beta}(\tilde{z}_0) - (1-\tilde{z}_0)\log(1-\tilde{z}_0) - \delta/2 \geq \alpha,$$

so

$$\liminf_{N \to \infty} \mathbb{P}\left( \frac{|V^{N*}|}{N} \geq 1 - \delta \text{ and } \frac{|\Lambda^{N*}|}{N} \geq \alpha \right) \geq \mathbb{P}(Z = 1)$$

as required.

**Acknowledgments.** The authors thank the Mathematisches Forschungsinstitut Oberwolfach for the invitation to the meeting in 2000 at which this collaboration began, and Brad Lackey for his interest in this project. David Levin helped us to articulate some of the concepts presented herein.

NATIONAL SECURITY AGENCY
P.O. BOX 535
ANNAPOLIS JUNCTION, MARYLAND 20701-0535
USA
E-MAIL: rwrd@afterlife.ncsc.mil

STATISTICAL LABORATORY
CENTRE FOR MATHEMATICAL SCIENCES
WILBERFORCE ROAD
CAMBRIDGE, CB3 0WB
UNITED KINGDOM
E-MAIL: j.r.norris@statslab.cam.ac.uk